\documentclass[11pt]{article}
\usepackage[active]{srcltx}

\usepackage{amsmath}
\usepackage{amsthm}
\usepackage{amsfonts}
\usepackage{amssymb}

\newcommand{\beq}{\begin{equation}}
\newcommand{\beqn}{\begin{equation*}}
\newcommand{\eeq}{\end{equation}}
\newcommand{\eeqn}{\end{equation*}}
\newcommand{\beqa}{\begin{eqnarray}}
\newcommand{\beqan}{\begin{eqnarray*}}
\newcommand{\eeqa}{\end{eqnarray}}
\newcommand{\eeqan}{\end{eqnarray*}}
\newcommand{\bdm}{\begin{displaymath}}
\newcommand{\edm}{\end{displaymath}}

\newcommand\nn{\nonumber}

\newcommand\benu{\begin{enumerate}}
\newcommand\eenu{\end{enumerate}}
\newcommand\bit{\begin{itemize}}
\newcommand\eit{\end{itemize}}

\def\End{\mathrm{End\,}}

\def\tr{\mathrm{tr\,}}

\def\dim{\mathrm{dim\,}}

\def\der'{\mathfrak{der}'\,}
\def\der{\mathfrak{der}\,}
\def\str'{\mathfrak{str}'\,}
\def\str{\mathfrak{str}\,}
\def\con{\mathfrak{con}\,}
\def\R{\mathbb{R}}

\def\C{\mathbb{C}}

\def\H{\mathbb{H}}

\def\bbO{\mathbb{O}}

\def\frake{\mathfrak{e}}

\def\g{\mathfrak{g}}
\def\h{\mathfrak{h}}
\def\so{\mathfrak{so}}

\def\e,n{e_1,\,e_2,\,\ldots,\,e_n}

\def\a,n{a_1,\,a_2,\,\ldots,\,a_n}

\def\al,n{\alpha_1,\,\alpha_2,\,\ldots,\,\alpha_n}

\newtheorem{theorem}{Theorem}

\theoremstyle{definition}

\begin{document}

\begin{center}

{\Large{\bf A generalization of the}}\\
\vspace{2mm}
{\Large{\bf Kantor-Koecher-Tits construction}}

\vspace{5mm}

Jakob Palmkvist

\footnotesize
\vspace{5mm}

{\it  Institut des Hautes Etudes Scientifiques\\ 35, Route de Chartres, FR-91440 Bures-sur-Yvette, 
France}

\vspace{10mm}

\it Talk presented at the Baltic-Nordic Workshop\\
\rm Algebra, Geometry, and Mathematical Physics,\\
\it G\"{o}teborg, \hbox{Sweden}, October 11--13, 2007.\\
Published in J.~Gen.~Lie Theory Appl.~{\bf 2}, (2008) 226.

\vspace{10mm}

\hrule

\vspace{5mm}

\parbox{130mm}{

\noindent \footnotesize \rm The Kantor-Koecher-Tits construction associates a Lie algebra to any Jordan algebra. We generalize this construction to include also extensions of the associated Lie algebra. In particular, the conformal realization of $\so(p+1,\,q+1)$ generalizes to $\so(p+n,\,q+n)$, for arbitrary $n$, with a linearly realized subalgebra $\so(p,\,q)$. We also show that the construction applied to $3 \times 3$ matrices over the division algebras $\R,\,\C,\,\H,\,\bbO$ gives rise to the exceptional Lie algebras $\mathfrak{f}_{4},\,\frake_{6},\,\frake_{7},\,\frake_{8}$, as well as to their affine, 
hyperbolic and further extensions.}

\vspace{5mm}

\hrule

\end{center}

\setcounter{equation}{0}
\setcounter{page}{1}

\noindent
This article aims to give a brief overview of results that have already been shown by the author in 
\cite{Palmkvist:2007as}, 
where details and a more comprehensive list of references are provided. We will mostly consider algebras over the real numbers, even though we will complexify real Lie algebras in order to study properties of the corresponding Dynkin diagrams. However, algebras are assumed to be real if nothing else is stated.

A {Jordan algebra} is a commutative algebra that satisfies the {Jordan identity}
\begin{align}
a^2 \circ (b \circ a) &= (a^2 \circ b) \circ a \label{the Jordan identity}
\end{align}
The symmetric part of the product in an associative but noncommutative algebra,
\begin{align*}
2 (a \circ b) =(ab+ba) 
\end{align*}
leads to a Jordan algebra in the same way as the antisymmetric part leads to a Lie algebra.
A deeper relationship between these two important kinds of algebras emerges in the generalization of a Jordan algebra to a {Jordan triple system} $J$
with a triple product 
$J^3 \rightarrow J$, $(x,\,y,\,z) \mapsto (xyz)$
such that $(xyz)=(zyx)$ and
\begin{align}
(uv(xyz))-(xy(uvz))=((uvx)yz)-(x(vuy)z)
\label{jtsidentitet}
\end{align}
This is indeed a generalization since any Jordan algebra, with triple product 
\begin{align*}
(xyz)=(x \circ y) \circ  z+x \circ (y \circ z)-y \circ (x \circ z)
\end{align*}
satisfies the conditions for a Jordan triple system. The same structure arises also in 
a 3-graded Lie algebra $\g=\g_{-1}+\g_0+\g_{1}$,
where $[\g_{m},\,\g_{n}] \subset \g_{m+n}$ (and $\g_k=0$ for $|k| > 1$). Let $\tau$ be a {graded} involution on $\g$, which means that $\tau(\g_k) \subset \g_{-k}$. Then the subspace 
$\g_{-1}$ is a Jordan triple system with the triple product
$(xyz)=[[x,\,\tau(y)],\,z]$.
Conversely, any Jordan triple system gives rise to a 3-graded Lie algebra $\g_{-1}+\g_0+\g_1$, spanned by the operators
\begin{align}
u &\in {\g_{-1}:} &   x&\mapsto u & &{{\text{(constant)}}}\nonumber\\
[u,\,\tau(v)]&\in{\g_{0}:} &   x&\mapsto(uvx) & &{\text{(linear)}}\nonumber\\
\tau(u) &\in {\g_{1}:} &   x&\mapsto-\tfrac{1}{2}(xux) & &{\text{(quadratic)}}
\label{jordanoperatorer}
\end{align}
where $x$ is an element in the Jordan triple system, which can be identified with $\g_{-1}$.
If the triple system is derived from a Jordan algebra $J$ by
(\ref{jtsidentitet}),
then the associated Lie algebra $\g_{-1}+\g_0+\g_1$ is 
the {conformal algebra}
$\con \, J$, and $\g_0$ is the {structure algebra} $\str \, J$.
If $J$ has an identity element, then all scalar multiplications form a one-dimensional ideal of $\str\, J$. Factoring out this ideal, we obtain the 
{reduced} structure algebra $\str' J$. Finally, all derivations of $J$ form a subalgebra $\der \, J$ of $\str' J$.
Thus,
$\g_0 = \str \, J \supset \str' J \supset \der \, J$. 
This construction of a 3-graded Lie algebra from a given Jordan algebra is called the {Kantor-Koecher-Tits construction} \cite{Kantor1,Koecher,Tits1}.
To see why the resulting Lie algebra is called 'conformal' we consider  
the Jordan algebras $H_2(\mathbb{K})$ of hermitian $2 \times 2$
matrices over the division algebras $\mathbb{K}=\R,\,\C,\,\mathbb{H},\bbO$. We have
\begin{align*} 
\nonumber
\der \,  H_2(\mathbb{K})&= \so(d-1)\\ \nonumber
\str'\,  H_2(\mathbb{K})&= \so(1,\,d-1)\\
\con \,  H_2(\mathbb{K})&= \so(2,\,d)
\end{align*}
for $d=3,\,4,\,6,\,10$, respectively [8]. It is well known that $\con \,  H_2(\mathbb{K})$ is the algebra that generates conformal transformations in a $d$-dimensional Minkowski spacetime. Furthermore, $\str'\,  H_2(\mathbb{K})$ is the Lorentz algebra and $\der \,  H_2(\mathbb{K})$ its spatial part. 
With a $d$-dimensional Minkowski spacetime we mean a vector space with a basis $P_{\mu}$ for $\mu=0,\,1,\,\ldots,\,d-1$ and an inner product such that $(P_{\mu},\,P_{\nu})=\eta_{\mu \nu}$, where $\eta=\text{diag}(-1,\,1,\,\ldots,\,1)$. We let a vector $x$ have components $x^{\mu}$ in this basis, $x=x^{\mu}P_{\mu}$, and let $\partial_{\mu}$ be the corresponding partial derivative. Then we can identify $\g_{-1}$ with this vector space and the operators  
(\ref{jordanoperatorer}) with the following vector fields:  
\begin{align}  
\nonumber
{P_{\mu}} &= {\partial_{\mu}} &&{\text{(translations)}}  \\
{G^{\mu}}_{\nu} &=
{x_{\nu}}{\partial^{\mu}}-{x^{\mu}}{\partial_{\nu}} &&{\text{(Lorentz transformations)}}\nonumber \\ \nonumber
{D} &= {x^{\mu}}{\partial_{\mu}} &&{\text{(dilatations)}}  \\
{K^{\mu}} &= -2{x^{\nu}} {x^{\mu}} {\partial_{\nu}} + 
{x^{\nu}} {x_{\nu}} {\partial^{\mu}} &&{\text{(special conformal
transformations)}}
\label{konfrel}
\end{align}
where the indices are raised and lowered by $\eta$. The fact that these operators satisfy the commutation relations for $\so(2,\,d)$ does not depend on the signature $(1,\,d-1)$ of $\eta$, so we have the same conformal realization of $\so(p+1,\,q+1)$ for any signature $(p,\,q)$. 

The Kantor-Koecher-Tits construction can be applied also to the Jordan algebras $H_3(\mathbb{K})$ of hermitian $3 \times 3$ matrices over $\mathbb{K}$
and then we obtain the first three rows in the 'magic square', which is a symmetric $4 \times 4$ array
of complex Lie algebras [8] (see also references therein). 
In particular, the exceptional Jordan algebra $H_3(\bbO)$ gives rise to exceptional Lie algebras:
\begin{align}
\der\, H_3(\bbO)&=\mathfrak{f}_4\nonumber 
\\ \nonumber
\str'\,             H_3(\bbO)&=\frake_{6}\\ 
\con\, H_3(\bbO)&=\frake_{7} \label{kkt}
\end{align}
(For simplicity, we do not specify the real forms of these complex Lie algebras.)
We focus on the $3 \times 3$ subarray in the lower right corner of the magic square, consisting of simply laced algebras:

\newcommand{\esex}{
\begin{picture}(70,30)(-5,-5)
\thicklines
\multiput(0,0)(15,0){5}{\circle{5}}
\multiput(2.5,0)(15,0){4}{\line(1,0){10}}
\put(30,2.5){\line(0,1){10}}
\put(30,15){\circle{5}}
\put(30,0){\circle*{5}}
\end{picture}}

\newcommand{\esju}{
\begin{picture}(85,30)(-5,-5)
\thicklines
\multiput(0,0)(15,0){6}{\circle{5}}
\multiput(2.5,0)(15,0){5}{\line(1,0){10}}
\put(45,2.5){\line(0,1){10}}
\put(45,15){\circle{5}}
\put(60,0){\circle*{5}}
\end{picture}}

\newcommand{\eatta}{
\begin{picture}(100,30)(-5,-5)
\thicklines
\multiput(0,0)(15,0){7}{\circle{5}}
\multiput(2.5,0)(15,0){6}{\line(1,0){10}}
\put(60,2.5){\line(0,1){10}}
\put(60,15){\circle{5}}
\put(15,0){\circle*{5}}
\end{picture}}

\newcommand{\esexminus}{
\begin{picture}(70,30)(-5,-5)
\thicklines
\multiput(0,0)(15,0){5}{\circle{5}}
\multiput(2.5,0)(15,0){4}{\line(1,0){10}}
\put(30,0){\circle*{5}}
\end{picture}}

\newcommand{\esjuminus}{
\begin{picture}(85,30)(-5,-5)
\thicklines
\multiput(0,0)(15,0){5}{\circle{5}}
\multiput(2.5,0)(15,0){4}{\line(1,0){10}}
\put(45,2.5){\line(0,1){10}}
\put(45,15){\circle{5}}
\put(60,0){\circle*{5}}
\end{picture}}

\newcommand{\eattaminus}{
\begin{picture}(100,30)(-5,-5)
\thicklines
\multiput(15,0)(15,0){6}{\circle{5}}
\multiput(17.5,0)(15,0){5}{\line(1,0){10}}
\put(60,2.5){\line(0,1){10}}
\put(60,15){\circle{5}}
\put(15,0){\circle*{5}}
\end{picture}}

\newcommand{\esexmminus}{
\begin{picture}(70,30)(-5,-5)
\thicklines
\multiput(0,0)(15,0){2}{\circle{5}}
\multiput(45,0)(15,0){2}{\circle{5}}
\multiput(2.5,0)(45,0){2}{\line(1,0){10}}
\end{picture}}

\newcommand{\esjumminus}{
\begin{picture}(85,30)(-5,-5)
\thicklines
\multiput(0,0)(15,0){4}{\circle{5}}
\multiput(2.5,0)(15,0){3}{\line(1,0){10}}
\put(45,2.5){\line(0,1){10}}
\put(45,15){\circle{5}}
\end{picture}}

\newcommand{\eattamminus}{
\begin{picture}(100,30)(-5,-5)
\thicklines
\multiput(30,0)(15,0){5}{\circle{5}}
\multiput(32.5,0)(15,0){4}{\line(1,0){10}}
\put(60,2.5){\line(0,1){10}}
\put(60,15){\circle{5}}
\end{picture}}

\newcommand{\bildC}{
\begin{picture}(30,30)(-5,-5)
\put(5,0){$\C$}
\end{picture}}

\newcommand{\bildH}{
\begin{picture}(30,30)(-5,-5)
\put(5,0){$\H$}
\end{picture}}

\newcommand{\bildbbO}{
\begin{picture}(30,30)(-5,-5)
\put(5,0){$\bbO$}
\end{picture}}

\newcommand{\bildKK}{
\begin{picture}(60,30)(-5,-5)
\put(20,0){$\mathbb{K}$}
\end{picture}}

\newcommand{\bildstr}{
\begin{picture}(60,30)(-5,-5)
\put(1,0){$\str'\,  H_3(\mathbb{K})$}
\end{picture}}

\newcommand{\bildcon}{
\begin{picture}(60,30)(-5,-5)
\put(1,0){$\con\,  H_3(\mathbb{K})$}
\end{picture}}

\begin{displaymath}
\qquad
\begin{array}{|c|c|c|c|}
\hline
\bildKK & \bildC & \bildH & \bildbbO \\
\hline 
\bildstr &\esexmminus & \esjumminus & \eattamminus \\ \hline
\bildcon &\esexminus & \esjuminus & \eattaminus \\ \hline
 & \esex & \esju & \eatta \\ \hline
\end{array}
\end{displaymath}\\
\label{kvadrat-sida}
{\!}It is easily seen that any simple root $\alpha$ of a complex Lie algebra $\g$ (or the corresponding node in the Dynkin diagram)
'generates' a grading of $\g$, where the subspace $\g_k$ is spanned by all root vectors $e_{\mu}$ or $f_{-\mu}$, such that the root $\mu$ has the coefficient 
$-k$ for $\alpha$ in the basis of simple roots [4].
In the middle row above,
the black node generates the 3-grading of the conformal algebra. (Here and below, this meaning of a black node in a Dynkin diagram should not be confused with any different meaning used elsewhere.)
The outermost node next to it in the last row generates
the unique 5-grading where 
the subspaces $\g_{\pm 2}$ are one-dimensional.
With this 5-grading, the algebras in the last row are called 'quasiconformal', associated to {Freudenthal triple systems} [2]. This is usually the way 
$\frake_8$ is included in the context of Jordan algebras and octonions.
The approach in this contribution is different: we want to generalize the conformal realization but keep the linear realization of the reduced structure algebra, and therefore we consider in the last row the grading generated by the black node itself. This grading seems better suited for a further generalization
to extensions of these exceptional Lie algebras.
We thus consider
the case when a finite Kac-Moody algebra 
$\h$ is extended to another one $\g$ in the following way, for an arbitrary integer $n \geq 2$.

\noindent
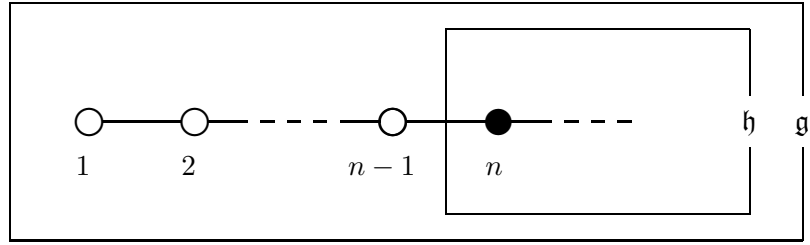
\begin{minipage}{260pt}
\begin{center}
\begin{picture}(190,80)(65,-10)
\put(125,-10){$1$}
\put(165,-10){$2$}
\put(228,-10){${n-1}$}
\put(280,-10){${n}$}
\put(265,10){\line(0,1){35}}
\put(265,10){\line(0,-1){35}}
\put(380,20){\line(0,1){25}}
\put(377,7){$\h$}
\put(380,0){\line(0,-1){25}}
\put(400,20){\line(0,1){35}}
\put(397,7){$\g$}
\put(400,0){\line(0,-1){35}}
\put(100,10){\line(0,1){45}}
\put(100,10){\line(0,-1){45}}
\put(100,55){\line(1,0){300}}
\put(100,-35){\line(1,0){300}}
\put(265,45){\line(1,0){115}}
\put(265,-25){\line(1,0){115}}
\thicklines
\multiput(130,10)(40,0){2}{\circle{10}}
\multiput(245,10)(40,0){1}{\circle{10}}
\multiput(285,10)(40,0){1}{\circle*{10}}
\multiput(135,10)(40,0){1}{\line(1,0){30}}
\multiput(250,10)(40,0){1}{\line(1,0){30}}
\put(175,10){\line(1,0){15}}
\put(230,10){\line(1,0){10}}
\multiput(195,10)(10,0){4}{\line(1,0){5}}
\put(290,10){\line(1,0){15}}
\multiput(310,10)(10,0){3}{\line(1,0){5}}
\put(245,10){\circle{10}}
\end{picture}
\end{center}
\vspace*{1.0cm}
\end{minipage} 

\noindent
The black node, which $\g$ and $\h$ have in common, generates a grading of $\g$ as well as of $\h$. 
We want to investigate how the triple systems $\g_{-1}$ and $\h_{-1}$, corresponding to these two gradings, are related to each other. 
It is clear that $\dim \g_{-1} = n\, \dim \h_{-1}$, which means that $\g_{-1}$ as a vector space is isomorphic to the direct sum $(\h_{-1})^n$ of $n$ vector spaces, each isomorphic to $\h_{-1}$. The question is if we can define a triple product on $(\h_{-1})^n$ such that $\g_{-1}$ and $(\h_{-1})^n$ are isomorphic also as triple systems. To answer this question, we 
write a general element in $(\h_{-1})^{n}$ as $(x_1)^1 + (x_2)^2 + \cdots + (x_{n})^{n}$, where $x_1,\,x_2,\,\ldots$ are elements in $\h_{-1}$.
Furthermore, we define for any graded involution $\tau$ on $\h$ a bilinear form on $\h_{-1}$ 
{associated to} $\tau$
by $(e_{\mu},\,\tau(f_{\nu}))=\delta_{\mu \nu}$ for root vectors
$e_{\mu} \in \h_{-1}$ and $f_{\nu} \in \h_1$.
The answer is then given by the following theorem 
(for a proof, see \cite{Palmkvist:2007as}).
\begin{theorem} \label{sats1}
The vector space $(\h_{-1})^{n}$, together with the triple product given by
\begin{align*}
(x^a y^b z^c)=
\delta^{ab}[[x,\,\tau(y)],\,z]^c  
- \delta^{ab}(x,\,y)z^c
+ \delta^{bc}(x,\,y)z^a
\end{align*}
for $a,\,b,\,\ldots=1,\,2,\,\ldots,\,n$ and $x,\,y,\,z \in \h_{-1}$,
is a triple system isomorphic to the triple system $\g_{-1}$ with the triple product 
$(uvw)= [[u,\,\tau(v)],\,w]$,
where the involution $\tau$ is extended from $\mathfrak{h}$ to $\g$ by $\tau(e^i)=-f^i$ for the simple root vectors.
\end{theorem}

\noindent
Even though the grading of $\h$ generated by the black node is a 3-grading, the grading of the extended algebra $\g$ generated by the same black node is an $m$-grading, where $m$ can be any odd positive integer, or even infinity.
Equivalently, even though the triple system $\h_{-1}$ is a Jordan triple system, this is in general not the case for the triple system $\g_{-1}$. However, $\g_{-1}$ is always a {generalized} Jordan triple system, which means that (\ref{jtsidentitet}) holds, but the triple product $(xyz)$ does not need to be symmetric in $x$ and $z$. The construction of the associated Lie algebra can be extended to any generalized Jordan triple system, as we will see an example of next.

When $\h=\so(p,\,q)$ for any positive integers $p,\,q$, the black node always generates a 5-grading of 
$\g=\so(p+n,\,q+n)$. Equivalently, $\g_{-1}$ is a generalized Jordan triple system {of second order}, or a {Kantor triple system} [1,\,4]. We get back the associated Lie algebra $\so(p+n,\,q+n)$ as the one spanned by the operators
\begin{align} 
[u,\,v]&\in \g_{-2}:&& z + Z\mapsto \langle u,\, v \rangle\nn\\
u&\in \g_{-1}:&& z + Z \mapsto u + \tfrac{1}{2} \langle u,\, z \rangle\nn\\
[u,\,\tau(v)]&\in \g_{0}:&& z + Z\mapsto (uvz) 
-\langle u,\, Z (v) \rangle\nn\\
\tau({u})&\in \g_{1}:&& z + Z \mapsto 
-\tfrac{1}{2}(zuz)- Z (u)\nn\\
                &&&\quad +\tfrac{1}{12}\langle (zuz),\,z \rangle-\tfrac{1}{2}\langle Z (u),\,z \rangle\nn\\
[\tau({u}),\,\tau(v)]&\in \g_{2}:&& z+Z \mapsto -\tfrac{1}{6}(z\langle u,\,v\rangle(z)z)
-Z(\langle u,\,v \rangle (z))\nn\\
                &&&\quad + \tfrac{1}{24} \langle (z\langle u,\,v\rangle(z)z),\, z \rangle
+\langle Z(u),\,Z(v) \rangle \label{kantoroperatorer}
\end{align}
where $z$ is an element in the Kantor triple system $\g_{-1}$, while $Z$ and the $\langle \,\ ,\ \rangle$ expressions belong to a certain subspace of $\End \g_{-1}$, which can be identified with $\g_{-2}$ \cite{Palmkvist:2005gc}. In the same way as (\ref{jordanoperatorer}), we can write the operators (\ref{kantoroperatorer}) as the following vector fields:
\begin{align*} \nonumber
P^{ab} &= -2{\partial}^{ab}\\ \nonumber
{P_{\mu}}^a &= {\partial_{\mu}}^a-2{x_{\mu b}}{\partial}^{ab}\\ \nonumber
{G^{\mu}}_{\nu} &= {x_{\nu a}}{\partial^{\mu a}}-{x^{\mu}}_a{\partial_{\nu}}^a\\ \nonumber
{D^a}_b &= {x^{\mu}}_b{\partial_{\mu}}^a+2x_{bc}\partial^{ac}\\ \nonumber
{K^{\mu}}_a &= -2{x^{\nu}}_a {x^{\mu}}_b {\partial_{\nu}}^b + 
{x^{\nu}}_a {x_{\nu b}} {\partial^{\mu b}}
-x_{ab} {\partial^{\mu b}}\\ \nonumber
            & \quad -2{x^{\nu}}_a {x^{\mu}}_b {x_{\nu c}} \partial^{bc} + 2x_{ab}{x^{\mu}}_c \partial^{bc}\\ \nonumber
K_{ab} &= {x^{\mu}}_a {x^{\nu}}_b {x_{\mu c}} {\partial_{\nu}}^c
-{x^{\mu}}_b {x^{\nu}}_a {x_{\mu c}} {\partial_{\nu}}^c - x_{ac} {x^{\mu}}_b {\partial_{\mu}}^c+ x_{bc} {x^{\mu}}_a {\partial_{\mu}}^c\\
            & \quad +
2{x^{\mu}}_a {x^{\nu}}_b {x_{\mu c}} {x_{\nu d}} \partial^{cd} - 2 x_{ac}x_{bd}\partial^{cd}
\end{align*}
When $(p,\,q)=(1,\,2),\,(1,\,3),\,(1,\,5),\,(1,\,9)$, we can 
thus construct $\so(p+n,\,q+n)$
starting from the Jordan algebra $H_2(\mathbb{K})$ and using the theorem.
In turns out that the bilinear form associated to the graded involution in this case is given by the trace: $(x,\,y)=\tr{(x \circ y)}$.
Then $H_2(\mathbb{K})^n$ will be a Kantor triple system with the triple product 
\begin{align} 
\label{tripppel}
(x^a y^b z^c)
&=
2\delta^{ab} ((z\circ y)\circ x)^c-2\delta^{ab}((z\circ x) \circ y )^c 
+2\delta^{ab}((x \circ y) \circ z)^c\nn\\
&\quad -\delta^{ab}(x,\,y)z^c+\delta^{bc}(x,\,y)z^a
\end{align}
where $a,\,b,\,c=1,\,2,\,\ldots,\,n$, and the Lie algebra associated to this Kantor triple system is thus $\so(p+n,\,q+n)$.
When we apply the same idea to the Jordan algebras $H_3(\mathbb{K})$, we find that
the exceptional algebras $\mathfrak{f}_4,\,\frake_6,\,\frake_7,\,\frake_8$
are the Lie algebras associated to $H_3(\mathbb{K})^2$, with the triple product
(\ref{tripppel})
for $\mathbb{K}=\R,\,\C,\,\H,\,\bbO$, respectively. Their affine and hyperbolic extensions are those associated to $H_3(\mathbb{K})^3$ and $H_3(\mathbb{K})^4$, respectively, while further extensions correspond to $H_3(\mathbb{K})^n$ for $n=5,\,6,\ldots$.
It remains to show that the bilinear form associated to the graded involution really is given by $(x,\,y)=\tr{(x \circ y)}$, also in the case of $3 \times 3$ matrices. However, 
one can show
that the triple product (\ref{tripppel}) indeed satisfies the definition of a 
generalized Jordan triple system when $x,\,y,\,z$ are elements in
$H_3(\mathbb{K})$
and $(x,\,y)=\tr{(x \circ y)}$.

An important and interesting difference between the $H_2(\mathbb{K})^n$ and
$H_3(\mathbb{K})^n$ cases is that the Lie algebra associated to $H_2(\mathbb{K})^n$ is 3-graded for $n=1$ and then 5-graded for all $n \geq 2$, while the Lie algebra associated to $H_3(\mathbb{K})^n$ is 3-graded for $n=1$ 
but 7-graded for $n = 2$,
and for $n=3,\,4,\,5,\ldots$,
we get infinitely many subspaces in the grading, since these Lie algebras are infinite-dimensional.
In the affine case, we only get the corresponding current algebra directly in this construction, which means that the central element and the derivation must be added by hand. It would be interesting to find an interpretation of these elements in the Jordan algebra approach. 
Finally, concerning the hyperbolic case and further extensions, we hope that our new construction can give more information about these indefinite Kac-Moody algebras, which, in spite of a great interest from both mathematicians and physicists, are not yet fully understood.

\subsection*{Acknowledgments} 
I would like to thank  the organizers of the AGMF workshop for the opportunity to present 
a talk at this nice conference. I am also grateful to Daniel Persson for comments on the manuscript.

\raggedright

\smallskip

\label{lastpage}

\end{document}